\documentclass[10pt,a4paper]{article}
\usepackage{amsfonts}
\usepackage{amsthm}
\usepackage{amssymb}
\usepackage[centertags]{amsmath}
\usepackage{amssymb}
\usepackage{enumerate}
\usepackage{graphicx}
\providecommand{\U}[2]{\protect\rule{1.5in}{1.5in}}
\theoremstyle{plain}
\newtheorem{theorem}{Theorem}[section]

\newtheorem{definition}[theorem]{Definition}
\newtheorem{example}[theorem]{Example}
\newtheorem{lemma}[theorem]{Lemma}

\numberwithin{equation}{section}
\begin{document}
\bibliographystyle{unsrt}
\title{Solving fuzzy convolution Volterra integro-differential equation via fuzzy Laplace transforms}
\author{Saif Ullah\footnote{Department of Mathematics, University of Peshawar, 25120, Khyber Pakhtunkhwa, Pakistan. E-mail: saifullah.maths@upesh.edu.pk}, Latif Ahmad\footnote{a. Shaheed Benazir Bhutto University, Sheringal, Dir, Khyber Pakhtunkhwa, Pakistan. b. Department of Mathematics, University of Peshawar, 25120, Khyber Pakhtunkhwa, Pakistan. E-mail: ahmad49960@yahoo.com}, Muhammad Farooq\footnote{Department of Mathematics, University of Peshawar, 25120, Khyber Pakhtunkhwa, Pakistan. E-mail: mfarooq@upesh.edu.pk}, Saleem Abdullah\footnote{Department of Mathematics, Quaid-i-Azam University, Islamabad, Pakistan. E-mail: saleemabdullah81@yahoo.com}}
\maketitle
\begin{abstract}
The solution of fuzzy integro-differential equations have a major role in the fields of science and engineering. Different approaches both numerical and analytic are used to solve these type of equations. In this paper, the solution of fuzzy convolution Volterra integro-differential equation  is investigated using fuzzy Laplace transform method under generalized Hukuhara differentiability. Finally, the method is illustrated with few examples to show the ability of the proposed method.
\end{abstract}
\section{Introduction}
One of the most important part of the fuzzy theory is the  fuzzy differential equations \cite{1:lt, 2:lt}, fuzzy integral equations \cite{3:lt, 4:lt, 5:lt} and fuzzy integro-differential equations (FIDEs) \cite{6:lt, 7:lt, 8:lt}. The  FIDEs is obtained when a physical system is modeled under differential sense \cite{26:lt}. Also FIDEs in fuzzy setting are a natural way to model uncertainty of dynamical systems. Therefore the solution of the fuzzy integro-differential equations is very important in various fields such as Physics, Geographic, Medical and Biological Sciences \cite{27:lt, 28:lt, 29:lt}. In \cite{9:lt} Seikkala defined fuzzy  derivatives and then some generalization have been investigated in \cite{10:lt, 11:lt, 12:lt, 13:lt}. The concept of integration of fuzzy functions was first introduced by Dubois and Prade \cite{11:lt}. Alternative approaches were later studied in \cite{14:lt, 15:lt}.

Different approaches are used to solve fuzzy integro-differential equation numerically and analytically \cite{16:lt,17:lt}.  In \cite{16:lt} the fuzzy differential transform method (FDTM) was used to solve fuzzy integro-differential equation (FIDE) of first order. The FDTM is semi analytic method which gives the exact solution if it is in the term of series expansion of some known function \cite{16:lt}. In \cite{17:lt} a series solution of fuzzy  Volterra integro-differential equation was studied.

In \cite{18:lt} Allahveranloo and Salashour  proposed the idea of fuzzy Laplace transform method for solving first order fuzzy differential equations under generalized H-differentiability. Recently Salashour developed the technique of fuzzy Laplace transform method to solve fuzzy convolution Volterra integral equation (FCVIE) of the second kind in \cite{19:lt}. In \cite{19:lt} original FCVIE is converted to two crisp convolution integral equations in order to determine the lower and upper function of solution, using fuzzy convolution operator.

\medskip
 In this paper we propose a technique for solving fuzzy convolution Volttera integro-differential equation (FCVIDE) using fuzzy Laplace transform. The general form of (FCVIDE) is
 \begin{equation}\label{s1}
x'(t)=f(t)+\int_0^t k(t-s)x(s)ds,
\end{equation}
with initial condition $x(0,r)=(\underline{x}, \overline{x})$.\\

\noindent The paper is organized as follows:\\
In section 2, some basic definitions and results are stated which will be used throughout this paper.
In section 3, the fuzzy Laplace transform and fuzzy convolution are briefly recalled. In section 4, the  fuzzy Laplace transforms are applied to fuzzy convolution Volttera integro-differential equation to construct the general technique. Illustrative examples are also considered to show the ability of the proposed method in section 5, and the conclusion is drawn in section 6.

\section{Preliminaries}
In this section we will recall some basics definitions and theorems needed throughout the paper such as fuzzy number, fuzzy-valued function and the derivative of the fuzzy-valued functions \cite{20:lt, 21:lt}.
\begin{definition} A fuzzy number is defined  as the mapping such that $u:R\rightarrow[0,1]$, which satisfies the following four properties
\begin{enumerate}
\item $u$ is upper semi-continuous.
 \item $u$ is fuzzy  convex that is $u(\lambda  x+(1-\lambda)y) \geq \min{\{u(x), u(y)\}}.\; x, y\in R$ and $\lambda\in [0,1]$.
\item $u$ is normal that is $\exists$ $x_0\in R$, where $u(x_0)=1$.
\item $A=\{\overline{x \in \mathbb{R}: u(x)>0}\}$ is compact, where $\overline{A}$ is closure of $A$.
\end{enumerate}
\end{definition}

\begin{definition}
A fuzzy number in parametric form is given as an order pair of the form $u=(\underline{u}(r), \overline{u}(r))$, where $0\leq r\leq1$ satisfying the following conditions.
\begin{enumerate}
\item $\underline{u}(r)$ is a bounded left continuous increasing function in the interval $[0,1]$.
\item $\overline{u}(r)$ is a bounded left continuous decreasing function in the interval $[0,1]$.
\item $\underline{u}{(r)\leq\overline{u}(r)}$. 
\end{enumerate}
If $\underline{u}(r)=\overline{u}(r)=r$, then $r$ is called crisp number.
\end{definition}

\noindent Since each $y\in R$ can be regarded as a fuzzy number if
\begin{eqnarray*}\widetilde{y}(t)=\begin{cases}1, \;\;\; if \;\; y=t,\\  0, \;\;\; if \;\; y\neq t.\end{cases}\end{eqnarray*}
For arbitrary fuzzy numbers $u=(\underline{u}(\alpha), \overline{u}(\alpha))$ and $v=(\underline{v}(\alpha), \overline{v}(\alpha))$ and an arbitrary crisp number $j$, we define addition and scalar multiplication as:
\begin{enumerate}
\item $(\underline{u+v})(\alpha)=(\underline{u}(\alpha)+\underline{v}(\alpha))$.
\item $(\overline{u+v})(\alpha)=(\overline{u}(\alpha)+\overline{v}(\alpha))$.
\item $(j\underline{u})(\alpha)=j\underline{u}(\alpha)$, $(j\overline{u})(\alpha)=j\overline{u}(\alpha)$, \mbox{       }  $j\geq0$.
\item $(j\underline{u})(\alpha)=j\overline{u}(\alpha)\alpha, (j\overline{u})(\alpha)=j\underline{u}(\alpha)\alpha$, $j<0$.
\end{enumerate}
\begin{definition} (See \cite{18:lt,23:lt}) Let us suppose that x, y $\in E$, if $\exists$ $z\in E$ such that
$x=y+z$, then $z$ is called the H-difference of $x$ and $y$ and is given by $x\ominus y$.\end{definition}

\noindent The Housdorff distance between the fuzzy numbers \cite{6:lt,12:lt,18:lt,23:lt} defined by
\[d:E\times E\longrightarrow R^{+}\cup \{{0}\},\]
\[d(u,v)=\sup_{r\in[0,1]}\max\{|\underline{u}(r)-\underline{v}(r)|, |\overline{u}(r)-\overline{v}(r)|\},\] \noindent where $u=(\underline{u}(r), \overline{u}(r))$ and $v=(\underline{v}(r), \overline{v}(r))\subset R$.
\\\\
We know that if $d$ is a metric in $E$, then it will satisfy the following properties, introduced by Puri and Ralescu \cite{22:lt}:
\begin{enumerate}
\item $d(u+w,v+w)=d(u,v)$, $\forall$  u, v, w $\in$ E.

\item $(k \odot u, k \odot v)=|k|d(u, v)$, $\forall$ k $\in$ R, \mbox{  and  } u, v $\in$ E.

\item $d(u \oplus v, w \oplus e)\leq d(u,w)+d(v,e)$, $\forall$ u, v, w, e $\in$  E.
\end{enumerate}

\begin{theorem} (see Wu \cite{24:lt}) Let $f$ be a fuzzy-valued function on $[a,\infty)$ given in the parametric form as $(\underline{f}(x,r), \overline{f}(x,r))$ for any constant number $r\in[0,1]$. Here we assume that $\underline{f}(x,r)$ and $\overline{f}(x,r)$ are Riemann-Integrable on $[a,b]$ for every $b\geq a$. Also we assume that $\underline{M}(r)$ and $\overline{M}(r)$ are two positive functions, such that
$\int_a^b|\underline{f}(x,r)| dx \leq \underline{M}(r)$ and $\int_a^b |\overline{f}(x,r)| dx \leq \overline{M}(r)$
for every $b\geq a$, then $f(x)$ is improper fuzzy Riemann-integrable on $[{a}, \infty)$. Thus an improper integral will always be a fuzzy number. In short \[ \int_a^b f(x) dx = \bigg( \int_a^b\underline{f}(x,r) dx, \int_a^b \overline{f}(x,r) dx\bigg).\]
It is will known that Hukuhare differentiability for fuzzy function was introduced by Puri \& Ralescu in \cite{22:lt}. 
\end{theorem}
\begin{definition}(see \cite{23:lt}) Let $f:(a,b)\rightarrow E$, where $x_{0}\in (a,b)$, then we say that $f$ is strongly generalized differentiable at $x_0$ (Beds and Gal differentiability).
If $\exists$ an element $f'(x_0)\in E$ such that
\begin{enumerate}
 \item $\forall h>0$ sufficiently small $\exists$ $f(x_0+h)\ominus f(x_0)$, $f(x_0)\ominus f(x_0-h)$, then the following limits hold (in the metric $d$)\\
 $\lim_{h\rightarrow 0}\frac{f(x_0+h)\ominus f(x_0)}{h}=\lim_{h\rightarrow 0}\frac{f(x_0)\ominus f(x_0-h)}{h}=f'(x_0)$,

\noindent or
\item $\forall h>0$ sufficiently small, $\exists$ $f(x_0)\ominus f(x_0+h)$,  $f(x_0-h)\ominus f(x_0)$, then the following limits hold (in the metric $d$) \\$\lim_{h\rightarrow 0}\frac{f(x_0)\ominus f(x_0+h)}{-h}=\lim_{h\rightarrow 0}\frac{f(x_0-h)\ominus f(x_0)}{-h}=f'(x_0)$,

    \noindent or
 \item $\forall h>0$ sufficiently small $\exists$ $f(x_0+h)\ominus f(x_0)$,  $f(x_0-h)\ominus f(x_0)$ and the following limits hold (in metric $d$)\\
$\lim_{h\rightarrow 0}\frac{(x_0+h)\ominus f(x_0)}{h}=\lim_{h\rightarrow 0}\frac{f(x_0-h)\ominus f(x_0)}{-h}=f'(x_0)$,

\noindent or
\item $\forall h>0$ sufficiently small $\exists$ $f(x_0)\ominus f(x_0+h)$,  $f(x_0)\ominus f(x_0-h)$, then the following limits holds(in metric $d$)\\
$\lim_{h\rightarrow 0}\frac{f(x_0)\ominus f(x_0+h)}{-h}=\lim_{h\rightarrow 0}\frac{f(x_0-h)\ominus f(x_0)}{h}=f'(x_0)$.
\end{enumerate}
The denominators $h$ and $-h$ denote multiplication by $\frac{1}{h}$ $\frac{-1}{h}$ respectively.\end{definition}
\begin{theorem}(See Chalco and Reman-Flores \cite{25:lt}) Let $f:R\rightarrow E$ be a function denoted by $f(t)=(\underline{f}(t,r),\overline{f}(t,r))$ for each $r\in[0,1]$. Then
\begin{enumerate}
\item If $f$ is $(i)$-differentiable, then $\underline{f}(t,r)$ and $\overline{f}(t,r)$ are differentiable functions and $f'(t)=(\underline{f}'(t,r), \overline{f}'(t,r))$.
\item If $f$ is $(ii)$-differentiable, then $\underline{f}(t,r)$ and $\overline{f}(t,r)$ are differentiable functions and $f'(t)=(\overline{f}'(t,r), \underline{f}'(t,r))$.
\end{enumerate}
\end{theorem}
\begin{lemma}(See Bede and Gal\cite{23:lt})
Let $x_0\in R$, then the FDE $y'=f(x,y)$, $y(x_0)=y_0\in R$ and $f:R\times E\rightarrow E$ is supposed to be a continuous function and equivalent to one of the following integral equations.
\[y(x)=y_0+\int_{x_0}^x f(t, y(t))dt \;\;\; \forall  \;\;\; x\in [x_0, x_1],\]
\noindent or
\[y(0)=y^1(x)+(-1)\odot\int_{x_0}^x f(t,y(t))dt\;\;\; \forall \;\;\; x\in [x_0, x_1],\]
\noindent on some interval $(x_0, x_1)\subset R$ depending on the strongly generalized differentiability. Integral equivalency shows that if one solution satisfies the given equation, then the other will also satisfy.
\end{lemma}

\section{Fuzzy Laplace transform and fuzzy convolution}
In this section we state some definitions and theorems from \cite{18:lt,19:lt} which will be used in the next section. Let $f$ is a fuzzy-valued function and $p$ is a real parameter, then according to \cite{18:lt} FLT of the function $f$ is defined as follows:
\begin{definition}

The FLT of fuzzy-valued function is
\begin{equation*}\label{eq2}F(p)=L[f(t)]=\int_{0}^{\infty}e^{-pt}f(t)dt=\lim_{\tau\rightarrow\infty}\int_{0}^{\tau}e^{-pt}f(t)dt,\end{equation*}
\begin{equation*}F(p)=\bigg[\lim_{\tau\rightarrow\infty}\int_{0}^{\tau}e^{-pt}\underline{f}(t)dt,\lim_{\tau\rightarrow\infty}\int_{0}^{\tau}e^{-pt}\overline{f}(t)dt\bigg],\end{equation*}
\noindent whenever the limits exist.
\end{definition}
\begin{equation*}F(p;r)=L[f(t;r)]=[l(\underline{f}(t;r)),l(\overline{f}(t;r))],\end{equation*}
\noindent where
\begin{equation*}l[\underline{f}(t;r)]=\int_{0}^{\infty}e^{-pt}\underline{f}(t;r)dt=\lim_{\tau\rightarrow\infty} \int_{0}^{\tau}e^{-pt}\underline{f}(t;r)dt,\end{equation*}
\begin{equation*}\label{eq6}l[\overline{f}(t;r)]=\int_{0}^{\infty}e^{-pt}\overline{f}(t;r)dt=\lim_{\tau\rightarrow\infty}\int_{0}^{\tau}e^{-pt}\overline{f}(t;r)dt.  \end{equation*}
\begin{definition}
The fuzzy convolution of two fuzzy-valued functions f and g defined by
\begin{equation*}
(f\ast g)(t)=\int_{0}^{t}f(s)g(t-s)ds,
\end{equation*}
where $t>0$ and it exists if f and g are say, piecewise continues functions.
\end{definition}
\begin{theorem}(Derivative theorem) Suppose that f is continues fuzzy valued function on $[0, \infty)$ and of exponential order $\alpha$ and that $f'$ is piecewise continues
 in $[0,\infty)$ then
\begin{enumerate}
\item \[L(f'(t))=pL(f(t))\circleddash f(0),\]
if f is (i)-differentiable.
\item \[L(f'(t))=(-f(0))\circleddash (-pL(f(t))),\]
if f is (ii)-differentiable.
\end{enumerate}
\end{theorem}
\begin{theorem}\label{ct}(Fuzzy convolution theorem),
if f and  g are piecewise continuous fuzzy-valued function on $[0,\infty)$, and of exponential order p, then
\begin{equation*}
L[(f\ast g)(t)]=L[f(t)].L[g(t)].
\end{equation*}
\end{theorem}
\section{Constructing the method}
 In this section, we will investigate solution of fuzzy convolution Volterra integro-differential equation using fuzzy Laplace transform. Consider the general equation (\ref{s1}), then applying fuzzy Laplace on both side  of (\ref{s1}) we get
\begin{equation}\label{s3}
L[x'(t)]=L[f(t)]+L\bigg[\int_0^t k(t-s)x(s)ds\bigg].\end{equation}
Using fuzzy  convolution theorem \ref{ct} on integral part we get
\begin{equation}\label{s33}
L[x'(t)]=L[f(t)]+L[k(t)]L[x(t)].
\end{equation}
\section*{Case (A1): When $x(t)$ is (i)-differentiable}
In this case from Theorem 2.6 and 3.3 we have $x'(t)=(\underline{x'}(t,r), \overline{x'}(t,r))$ and
\begin{equation}\label{s333}L(x'(t))=pL(x(t))\circleddash x(0).\end{equation}
\noindent Using (\ref{s333}) in (\ref{s33}), gives
\begin{equation}\label{s3333}pL(x(t))\circleddash x(0)=L[f(t)]+L[k(t)]L[x(t)].\end{equation}
Now the classical form of (\ref{s3333}) is given as follows:
\begin{equation}\label{s4}
pl[\underline{x}(t,r)]-\underline{x}(0,r)=l[\underline{f}(t,r)]+\underline{l[k(t,r)]l[x(t,r)]}.
\end{equation}
Similarly

\begin{equation}\label{s2}
pl[\overline{x}(t,r)]-\overline{x}(0,r)=l[\overline{f}(t,r)]+\overline{l[k(t,r)]l[x(t,r)]}.
\end{equation}
Now different cases arise
\begin{enumerate}
\item   When both $x(t,r)$ and $k(t,r)$ are positive then
\begin{eqnarray*}
\underline{l[k(t,r)]l[x(t,r)]}=l[\underline{k}(t,r)]l[\underline{x}(t,r)],\\
\overline{l[k(t,r)]l[x(t,r)]}=l[\overline{k}(t,r)]l[\overline{x}(t,r)].
\end{eqnarray*}
\item When both $x(t,r)$ and $k(t,r)$ are negative then
\begin{eqnarray*}
\underline{l[k(t,r)]l[x(t,r)]}=l[\overline{k}(t,r)]l[\overline{x}(t,r)],\\
\overline{l[k(t,r)]l[x(t,r)]}=l[\underline{k}(t,r)]l[\underline{x}(t,r)].
\end{eqnarray*}
\item When $k(t,r)$ is negative and $x(t,r)$ is positive then
\begin{eqnarray*}
\underline{l[k(t,r)]l[x(t,r)]}=l[\overline{k}(t,r)]l[\underline{x}(t,r)],\\
\overline{l[k(t,r)]l[x(t,r)]}=l[\underline{k}(t,r)]l[\overline{x}(t,r)].
\end{eqnarray*}
\item When $x(t,r)$ is negative and $k(t,r)$ is positive then
\begin{eqnarray*}
\underline{l[k(t,r)]l[x(t,r)]}=l[\underline{k}(t,r)]l[\overline{x}(t,r)],\\
\overline{l[k(t,r)]l[x(t,r)]}=l[\overline{k}(t,r)]l[\underline{x}(t,r)].
\end{eqnarray*}
\end{enumerate}
Here we discuss the first case in detail the other cases can be dealt in a similar way. When both $x(t,r)$ and $k(t,r)$ are positive then equation (\ref{s4}) becomes
\begin{equation}
pl[\underline{x}(t,r)]-\underline{x}(0,r)=l[\underline{f}(t,r)]+l[\underline{k}(t,r)]l[\underline{x}(t,r)],
\end{equation}
\begin{equation}
(p-l[\underline{k}(t,r)])l[\underline{x}(t,r)]=\underline{x}(0,r)+l[\underline{f}(t,r)],
\end{equation}
\begin{equation*}
l[\underline{x}(t,r)]=\frac{\underline{x}(0,r)+l[\underline{f}(t,r)]}{(p-l[\underline{k}(t,r)])}.
\end{equation*}

Similarly in upper case
\begin{equation*}
l[\overline{x}(t,r)]=\frac{\overline{x}(0,r)+l[\overline{f}(t,r)]}{(p-l[\overline{k}(t,r)])}.
\end{equation*}
Finally applying inverse Laplace to get required lower and upper solution as follows:
\begin{equation*}
\underline{x}(t,r)=l^{-1}\bigg[\frac{\underline{x}(0,r)+l[\underline{f}(t,r)]}{(p-l[\underline{k}(t,r)])}\bigg],
\end{equation*}
\begin{equation*}
\overline{x}(t,r)=l^{-1}\bigg[\frac{\overline{x}(0,r)+l[\overline{f}(t,r)]}{(p-l[\overline{k}(t,r)])}\bigg].
\end{equation*}
\section*{Case (A2): When $x(t)$ is (ii)-differentiable}
In this case using Theorem 2.6 and 3.3 we have $x'(t)=(\overline{x'}(t,r), \underline{x'}(t,r))$ and
\begin{equation}\label{s51}L(x'(t))=-x(0))\circleddash (-pL(x(t)).\end{equation}
\noindent Using (\ref{s51}) in (\ref{s33}), gives
\begin{equation}\label{s52}-x(0))\circleddash (-pL(x(t))=L[f(t)]+L[k(t)]L[x(t)].\end{equation}
 The classical form of (\ref{s52}) is as under:

 \begin{equation} \label{s7}
 pl[\underline{x}(t,r)]-\underline{x}(0,r)=l[\overline{f}(t,r)]+\overline{l[k(t,r)]l[x(t,r)]},\end{equation}
 \begin{equation}\label{s6}
 pl[\overline{x}(t,r)]-\overline{x}(0,r)=l[\underline{f}(t,r)]+\underline{l[k(t,r)]l[x(t,r)]}.\end{equation}
 Now for the case when both $ x(t,r)$ and $k(t,r)$ are positive then equations (\ref{s7}) and (\ref{s6}) can be written as follow:

 \begin{equation} \label{s9}
 pl[\underline{x}(t,r)]-\underline{x}(0,r)=l[\overline{f}(t,r)]+l[\overline{k}(t,r)]l[\overline{x}(t,r)],\end{equation}
 \begin{equation} \label{s8}
 pl[\overline{x}(t,r)]-\overline{x}(0,r)=l[\underline{f}(t,r)]+l[\underline{k}(t,r)]l[\underline{x}(t,r)].\end{equation}
 Now solving (\ref{s9}) and (\ref{s8}) for   $l[\underline{x}(t,r)]$ and $l[\overline{x}(t,r)]$ we get:

 \begin{eqnarray*}
 l[\underline{x}(t,r)]=\frac{l[\overline{k}(t,r)\overline{x}(0,r)]}{p^2-l[\overline{k}(t,r)]l[\underline{k}(t,r)]}+
 \frac{l[\overline{k}(t,r)]l[\underline{f}(t,r)]}{p^2-l[\overline{k}(t,r)]l[\underline{k}(t,r)]}\\
 +\frac{\underline{x}(0,r)p}{p^2-l[\overline{k}(t,r)]l[\underline{k}(t,r)]}+\frac{l[\overline{f}(t,r)]p}{p^2-l[\overline{k}(t,r)]l[\underline{k}(t,r)]}, \end{eqnarray*}
 \begin{eqnarray*}
 l[\overline{x}(t,r)]=\frac{l[\underline{k}(t,r)\underline{x}(0,r)]}{p^2-l[\underline{k}(t,r)]l[\overline{k}(t,r)]}+
 \frac{l[\underline{k}(t,r)]l[\overline{f}(t,r)]}{p^2-l[\underline{k}(t,r)]l[\overline{k}(t,r)]}\\
 +\frac{\overline{x}(0,r)p}{p^2-l[\underline{k}(t,r)]l[\overline{k}(t,r)]}+\frac{l[\underline{f}(t,r)]p}{p^2-l[\underline{k}(t,r)]l[\overline{k}(t,r)]}. \end{eqnarray*}
 Finally taking the inverse Laplace we get the lower and upper solutions of given integro-differential equation (\ref{s1}) which are:

 \begin{eqnarray*}
 \underline{x}(t,r)=l^{-1}\bigg[\frac{l[\overline{k}(t,r)\overline{x}(0,r)]}{p^2-l[\overline{k}(t,r)]l[\underline{k}(t,r)]}\bigg]+
 l^{-1}\bigg[\frac{l[\overline{k}(t,r)]l[\underline{f}(t,r)]}{p^2-l[\overline{k}(t,r)]l[\underline{k}(t,r)]}\bigg]\\
 +l^{-1}\bigg[\frac{\underline{x}(0,r)p}{p^2-l[\overline{k}(t,r)]l[\underline{k}(t,r)]}\bigg]+l^{-1}\bigg[\frac{l[\overline{f}(t,r)]p}{p^2-l[\overline{k}(t,r)]l[\underline{k}(t,r)]}\bigg], \end{eqnarray*}
 and
\begin{eqnarray*}
\overline{x}(t,r)=l^{-1}\bigg[\frac{l[\underline{k}(t,r)\underline{x}(0,r)]}{p^2-l[\underline{k}(t,r)]l[\overline{k}(t,r)]}\bigg]+
l^{-1}\bigg[\frac{l[\underline{k}(t,r)]l[\overline{f}(t,r)]}{p^2-l[\underline{k}(t,r)]l[\overline{k}(t,r)]}\bigg]\\
+l^{-1}\bigg[\frac{\overline{x}(0,r)p}{p^2-l[\underline{k}(t,r)]l[\overline{k}(t,r)]}\bigg]+l^{-1}\bigg[\frac{l[\underline{f}(t,r)]p}{p^2-l[\underline{k}(t,r)]l[\overline{k}(t,r)]}\bigg].
\end{eqnarray*}
\section{Numerical examples}
In this section we will discuss the solution of fuzzy convolution Volterra integro-differential equations using FLT to show the utility of the proposed method in Section 4.
\begin{example}
Let us consider the following fuzzy convolution Volterra integro-differential equation.
\begin{eqnarray}\label{l1}
x'(t)=(1+t)(r+1,r-2)+\int_0^{t}x(s)ds,
\\ \nonumber with \mbox{  } initial \mbox{  } condition \mbox{  } x(0,r)=(0,0).
\end{eqnarray}
\subsection{Case $A1$: (i)-differentiability}
Now if we apply FLT on (\ref{l1}), then we get
\begin{equation}\label{l2}
pL[x(t)]\ominus{x}(0)=(r+1,r-2)L[(1+t)]+L[k(t)]L[x(t)],
\end{equation}
But here in this case $k(t-s)=1$, then (\ref{l2}) becomes
\begin{equation}\label{l3}
pL[x(t)]\ominus{x}(0)=(r+1,r-2)L[(1+t)]+L[1]L[x(t)].
\end{equation}
Now the classical form of (\ref{l3}) is given by
\begin{equation}\label{l11}
pl[\underline{x}(t,r)]-\underline{x}(0,r)=(r+1)l(1+t)+\underline{l[1]l[x(t,r)]},
\end{equation}

and
\begin{equation}\label{l12}
pl[\overline{x}(t,r)]-\overline{x}(0,r)=(r-2)l(1+t)+\overline{l[1]l[x(t,r)]}.
\end{equation}
Now let $x(t,r)$ is positive then equations (\ref{l11}) and (\ref{l12}) become:
\begin{equation}\label{l13}
pl[\underline{x}(t,r)]-\underline{x}(0,r)=(r+1)l(1+t)+l[1]l[\underline{x}(t,r)],
\end{equation}
\begin{equation}\label{l14}
pl[\overline{x}(t,r)]-\overline{x}(0,r)=(r-2)l(1+t)+l[1]l[\overline{x}(t,r)].
\end{equation}
Which gives the final lower and upper solutions as given below:
\begin{equation*}
\underline{x}(t,r)=(r+1)[e^t-1],
\end{equation*}
\begin{equation*}
\overline{x}(t,r)=(r-2)[e^t-1].
\end{equation*}

Now if $x(t,r)$ is negative, then (\ref{l11}) and (\ref{l12}) become:

\begin{equation}\label{l4}
pl[\underline{x}(t,r)]-\underline{x}(0,r)=(r+1)l(1+t)+l[1]l[\overline{x}(t,r)],
\end{equation}
\begin{equation}\label{l5}
pl[\overline{x}(t,r)]-\overline{x}(0,r)=(r-2)l[(1+t)]+l[1]l[\underline{x}(t,r)].
\end{equation}
Solving (\ref{l4}) and (\ref{l5}) simultaneously, we get
\begin{equation*}
\underline{x}(t,r)=(r+1)\bigg[\frac{1}{2}(e^t+\sin(t)-\cos(t))\bigg]+(r-2)\bigg[\frac{1}{2}(e^t+\cos(t)-\sin(t))-1\bigg],
\end{equation*}
\begin{equation*}
\overline{x}(t,r)=(r-2)\bigg[\frac{1}{2}(e^t+\sin(t)-\cos(t))\bigg]+(r+1)\bigg[\frac{1}{2}(e^t+\cos(t)-\sin(t))-1\bigg].
\end{equation*}
\end{example}
\begin{example}
Let us consider the following FCVIDE
\begin{eqnarray}\label{l6}
x'(t)=(r-1,1-r)+\int_0^{t}x(s)ds,\\
 \nonumber with \mbox{  } initial \mbox{  } condition \mbox{  } x(0,r)=(0,0).
\end{eqnarray}
\subsection{Case $A2$: (ii)-differentiability}
Applying FLT on (\ref{l6}), we have
\begin{equation}
-(x(0))\ominus(-p)L[x(t)]=(r-1,1-r)L[1]+L\bigg[\int_0^tx(s)ds\bigg].
\end{equation}
The classical form by using characterization theorem as
\begin{equation}\label{l15}
pl[\underline{x}(t,r)]-\underline{x}(0,r)=(1-r)l[1]+\overline{l[1]l[{x}(t,r)]},
\end{equation}

\begin{equation}\label{l16}
pl[\overline{x}(t,r)]-\overline{x}(0,r)=(r-1)l(1)+\underline{l[1]l[x(t,r)]}.
\end{equation}

If $x(t,r)$ is negative then (\ref{l15}) and  (\ref{l16}) become
\begin{equation*}
pl[\underline{x}(t,r)]-\underline{x}(0,r)=(1-r)l(1)+l[1]l[\underline{x}(t,r)],
\end{equation*}
\begin{equation*}
pl[\overline{x}(t,r)]-\overline{x}(0,r)=(r-1)l(1)+l[1]l[\overline{x}(t,r)].
\end{equation*}
Solving for  $\underline{x}(t,r)$ and $\overline{x}(t,r)$ we get
\begin{equation*}
\underline{x}(t,r)=(1-r)\sinh(t),
\end{equation*}
\begin{equation*}
\overline{x}(t,r)=(r-1)\sinh(t).
\end{equation*}
Now we will discuss case when $x(t,r)$ is positive. Here equations (\ref{l15}) and  (\ref{l16}) simplify to
\begin{equation*}
pl[\underline{x}(t,r)]-\underline{x}(0,r)=(1-r)l(1)+l[1]l[\overline{x}(t,r)],
\end{equation*}
\begin{equation*}
pl[\overline{x}(t,r)]-\overline{x}(0,r)=(r-1)l(1)+l[1]l[\underline{x}(t,r)].
\end{equation*}
Solving  simultaneously for  $\underline{x}(t,r)$ and $\overline{x}(t,r)$, we get the following solutions:
\begin{equation*}
\underline{x}(t,r)=(1-r)\bigg[\frac{1}{2}\sinh(t)+\frac{1}{2}\sin(t)\bigg]+(r-1)\bigg[\frac{1}{2}\sinh(t)-\frac{1}{2}\sin(t)\bigg],
\end{equation*}
\begin{equation*}
\overline{x}(t,r)=(r-1)\bigg[\frac{1}{2}\sinh(t)+\frac{1}{2}\sin(t)\bigg]+(1-r)\bigg[\frac{1}{2}\sinh(t)-\frac{1}{2}\sin(t)\bigg].
\end{equation*}
\end{example}
\begin{example}
Let us consider the following FCVIDE.
\begin{eqnarray}\label{b2}
x'(t)=(r-1,1-r)+\int_0^{t}e^{-2x}x(t)dx,
\\ \nonumber with \mbox{  } initial \mbox{  } condition \mbox{  } u(0,r)=(r-1,1-r).
\end{eqnarray}
Now applying the FLT on both sides of (\ref{b2})
\begin{equation}
pL[x(t)]\ominus{x}(0)=(r-1,1-r)L[1]+L[e^{-2t}]L[x(t)].
\end{equation}
In case of (i)-differentiability (i.e Case $A1$) and when $x(t,r)$ is positive  the lower solution is given by
\begin{equation}
pl[\underline{x}(t,r)]-\underline{x}(0,r)=(r-1)l[1]+l[e^{-2t}]l[\underline{x}(t,r)],
\end{equation}
\begin{equation*}
\underline{x}(t,r)=(r-1)[e^{-t}+te^{-t}-3\delta(t)+2].
\end{equation*}
Also the upper solution is given by
\begin{equation}
pl[\overline{x}(t,r)]-\overline{x}(0,r)=(1-r)l[1]+l[e^{-2t}]l[\overline{x}(t,r)],
\end{equation}
\begin{equation*}
\overline{x}(t,r)=(1-r)[e^{-t}+te^{-t}-3\delta(t)+2].
\end{equation*}
\end{example}
\section{Conclusion}
The FLT method was applied for the solution of FCVIE in \cite{19:lt}. Here in this paper we investigated the applicability of fuzzy Laplace transform for the solution of FCVIDE. We have illustrated the method by solving some examples. This approach toward the solution of FCVIDE is simple and easy to understand. This work can possibly be extended to an $nth$ order FCVIDE. This work is in progress.

\bibliography{references}
\end{document}